\documentstyle[12pt]{amsart} \numberwithin{equation}{section}
\input{amssym.def}
\input{amssym.tex}

\allowdisplaybreaks
\begin{document}

\title[Hyperstability of orthogonally 3-Lie homomorphism] {Hyperstability of orthogonally 3-Lie homomorphism: an orthogonally fixed point approach}
\author[V. Keshavarz, S. Jahedi]
{Vahid Keshavarz$^*$, Sedigheh Jahedi}
\address{Vahid Keshavarz and Sedigheh Jahedi \newline   Department of Mathematics, Shiraz University of Technology, P. O. Box 71557-13875, Shiraz, Iran.}
\email{v.keshavarz68@@yahoo.com}
\email{jahedi@@sutech.ac.ir}
\author[Th. M. Rassias]
{Themistocles M. Rassias}
\address{Themistocles M. Rassias \newline   Department of Mathematics, National Technical University of Athens, Zografou Campus, 15780
Athens, Greece.}
\email{trassias@@math.ntua.gr}
\begin{abstract}
In this paper, by using the orthogonally fixed point method, we prove the Hyers-Ulam stability and the hyperstability of orthogonally 3-Lie homomorphisms for additive $\rho$-functional equation in 3-Lie algebras.\\
Indeed, we investigate the Hyers-Ulam stability and the hyperstability of the system of functional equations
\begin{eqnarray*}
\left\{
  \begin{array}{ll}
    f(x+y)-f(x)-f(y)= \rho(2f(\frac{x+y}{2})+ f(x)+ f(y)),\\
\\
    f([[x,y],z])=[[f(x),f(y)],f(z)]
  \end{array}
\right.
\end{eqnarray*}
in 3-Lie algebras
(where $\rho$ is a fixed real  number with $\rho \ne 1$).
\end{abstract}
\footnotetext{*Corresponding author}
\footnotetext {2010 Mathematics Subject Classification. Primary 46S10, 39B62, 39B52, 47H10, 12J25}
\footnotetext {Keywords: orthogonally fixed point method; Hyers-Ulam stability; 3-Lie homomorphism in orthogonally 3-Lie algebras; hyperstability; additive $\rho$-Jensen functional equation.}

\maketitle

\baselineskip=16pt

\theoremstyle{definition}
  \newtheorem{df}{Definition}[section]
  \newtheorem{remark}[df]{Remark}
\theoremstyle{plain}
  \newtheorem{lem}[df]{Lemma}
  \newtheorem{thm}[df]{Theorem}
  \newtheorem{pro}[df]{Proposition}
  \newtheorem{cor}[df]{Corollary}

\setcounter{section}{0}

\section{Introduction}
The stability problem of functional equations  originated from a
question of Ulam \cite{ul60}  concerning the stability of group
homomorphisms.
\\Hyers \cite{hy41} gave a first affirmative partial answer to the question of Ulam for Banach spaces.
Hyers' Theorem  was generalized by Th. M. Rassias \cite{ra78} for linear mappings by
considering an unbounded Cauchy difference.  A
generalization of the Rassias' Theorem was obtained by G\u
avruta \cite{ga94} by replacing the unbounded Cauchy difference by a
general control function in the spirit of  Rassias' approach.
\\In 1996, G. Isac and Th. M. Rassias \cite{ir} were the first to
provide applications of stability theory of functional equations for
the proof of new fixed point theorems with applications.
The stability problems of several functional
equations have been extensively investigated by a number of
authors (see \cite{24, ii2, plx, mt, mt1, mt2, mt3, mt4, mt5}).
\\There are several orthogonality notions on a real normed space
such as Birkhoff-James, Boussouis, (semi--)inner product, Singer,
Carlsson, unitary--Boussouis, Roberts, Pythagorean and Diminnie
(see, \cite{Alonso1, Alonso2}). But here, we present the
introduced orthogonally sets and some corresponding concepts by Eshaghi Gordji { et al.} in \cite{v3}. Recently , several authors worked on orthogonally fixed point(see \cite{m.r1, m.r, SDC}).\\
Eshaghi Gordji {et al.} in \cite{v3} introduced orthogonal sets and some corresponding concepts.
\begin{df}
$(i):$ Let $X\neq\emptyset$ and $\bot\subseteq X \times X$ be a binary
relation. If $\bot$ satisfies the following condition
\begin{equation*}
\exists x_0; (\forall y; y \bot x_{0})\hspace{.2cm} or\hspace{.2cm} (\forall y; x_{0}\bot y ),
\end{equation*}
it is called an orthogonal set (briefly O-set). We denote this O-set by $(X,\bot)$.\\
$(ii):$ Let $(X,\bot)$ be an O-set. A sequence $\{x_{n}\}_{n\in\Bbb{N}}$ is
called an orthogonal sequence (briefly O-sequence)  if
\begin{equation*}
(\forall n; x_{n}\bot x_{n+1} )\hspace{.2cm} or \hspace{.2cm} (\forall n; x_{n+1}\bot x_{n} ).
\end{equation*}
$(iii):$ If $(X,\bot)$ is an O-set and  $(X,d)$ is a metric space then $(X,\bot,d)$ is an orthogonally metric space. A mapping $f:X\to X$ is a $\bot-$continuous in $x\in X$ if for each O-sequence $\{x_{n}\}_{n\in \Bbb{N}}$ in $X$ with $x_{n}\to x$, $f(x_{n})\to f(x)$. Obviously, every continuous mapping is $\bot-$continuous.\\
$(iv):$ Similarly, a Cauchy sequence $\{x_n\}$ in $X$ is said to be a cauchy orthogonally sequence ( briefly, cauchy O-sequence) if $\forall n\in A,$ $x_n\bot x_{n+1}$ or $x_{n+1}\bot x_{n}.$ An orthogonal space $(X,\bot,d)$ is orthogonally complete (briefly O-complete) if every Cauchy O-sequence is convergent.

It is easy to see that every complete metric space is O-complete and
the converse is not true.

$(v):$ Let $(X,\bot,d)$ be an orthogonally metric space and $0< \lambda<1$.
A mapping  $f:X\to X$ is said to be orthogonality contraction
with Lipschitz constant $\lambda$ if
\begin{equation*}
 d(fx,fy)\leq \lambda d(x,y)\;\;\  if ~x \bot y.
\end{equation*}
\end{df}\vskip 2mm
Remember that a Lie algebra is a Banach  algebra  endowed with the Lie product
$$[x, y] := \frac{(xy- yx)}{2}.$$
Similarly, a 3-Lie algebra is a Banach algebra  endowed with the product
$$\Big[[x, y],z\Big] := \frac{[x, y]z- z[x, y]}{2}.$$
for all $x,y,z \in A$.
\begin{df}
Let $A$ and $B$ be two 3-Lie algebras. A  mapping $H: A \to B$ is called a orthogonally 3-Lie homomorphism if\\
(i) $H$ is a linear mapping.\\
(ii) For all $x,y,z \in A$ with $ x\bot y,~~ x \bot z,~~ y\bot z.$
 $$H([[x,y],z])=[[H(x),H(y)],H(z)].$$
\end{df}
Recently, Eshaghi Gordji {\it et al.} in \cite{v3} proved a fixed point theorem in O-sets as follows:
\begin{thm}\cite{v4}\label{T.1}
Let $(X,d,\bot)$ be an O-complete generalized metric space. Let $T:X\to X$ be
a $\perp$-preserving, $\perp$-continuous and $\bot$-$\lambda$-contraction.
Let $x_0 \in X$ satisfy for all $y\in X$, $x_0\perp y$
or for all $y\in X$, $y\perp x_0$, and consider the ``{\it O-sequence of successive approximations.
with initial element $x_0$}\,'': $x_0$, $T(x_0)$, $T^2(x_0)$, ..., $T^n(x_0)$, ... .
Then, either $d(T^n(x_0), T^{n+1}(x_0))= \infty$ for all
$n\geq 0$, or there exists a positive integer $n_0$ such that $d(T^n(x_0), T^{n+1}(x_0))<\infty$
for all $n>n_0$. If the second alternative holds, then\\
(i): the O-sequence of $\{T^n(x_0)\}$ is convergent to a fixed point $x^*$ of $T$\\
(ii): $x^*$ is the unique fixed point of $T$ in $X^*=\{y\in X : d(T^n(x_0), y) < \infty\}$\\
(iii): if $y \in X$, then
$$d(y, x^*) \leq{1\over 1- \lambda}d(y, T(y)).$$
\end{thm}\vskip 2mm
\section{Main results}
Throughout this Section, assume that $A$ and $B$ are two orthogonally 3-Lie algebras and let $\rho$ be a fixed real  number with $\rho \ne 1$ and for simplicity, denote

\begin{equation*}
\begin{split}
V_\bot:=\{x,y,z \in A\quad |\quad x\bot y,~~ x \bot z,~~ y\bot z\}
\end{split}
\end{equation*}
\begin{equation}\label{kv2}
\begin{split}
\Delta_\rho f(x,y):=f(x+y)-f(x) -f(y)  - \rho\left(2f\left(\frac{x+y}{2}\right) - f(x) - f(y)\right)
\end{split}
\end{equation}
\begin{equation}\label{kv3}
\begin{split}
\Re f([[x,y],z]):=f([[x,y],z]-[[f(x),f(y)],f(z)]
\end{split}
\end{equation}
where $x,y,z \in V_\bot .$
\begin{lem}\label{l2.1}\cite{psl}
Let $X$ and $Y$ be vector spaces.
If a mapping  $f : X \rightarrow Y$ satisfies
\begin{eqnarray}\label{2.1}
f(x+y)-f(x) -f(y)= \rho \left(2f\left(\frac{x+y}{2}\right)- f(x) - f(y)\right)
\end{eqnarray} for all $x, y \in X$, then $f : X \rightarrow Y$ is additive.
\end{lem}
In the following theorem, we prove the Hyers-Ulam stability of orthogonally 3-Lie homomorphism in orthogonally 3-Lie algebras.
\begin{thm}\label{t2.2}
Let $f:A\to B$ be a mapping and let $\varphi: A^2 \to [0,\infty)$ and $\psi: A^3 \to [0,\infty)$ be two functions such that there exists an $L<1$ with
\begin{eqnarray}\label{2.2}
\|\Delta_\rho f(x,y)\|\le \varphi(x,y)
\end{eqnarray}
and
\begin{eqnarray}\label{2.3}
\|\Re f([[a,b],c])\| \le \psi(a,b,c)
\end{eqnarray}
for all $x,y \in V_\bot$. If there exists a constant $0<L<1$ such that
\begin{eqnarray}\label{2.4}
 \varphi\left(\frac{x}{2}, \frac{y}{2}\right)
 \le \frac{L}{2} \varphi\left({x}, y\right)
\end{eqnarray}
\begin{eqnarray}\label{2.5}
 \psi\left(\frac{a}{2}, \frac{b}{2},\frac{c}{2}\right)
 \le \frac{L}{2^3} \psi\left(a,b,c\right)
\end{eqnarray}
 for all $x, y, a, b, c\in V_\bot$, then there exists a unique orthogonally 3-Lie homomorphism $\Im:A \rightarrow B$
such that
\begin{eqnarray}
\|f(x)-\Im(x) \|  \le  \frac{L}{2(1-L)} \varphi\left({x}, x\right)
\end{eqnarray} for all $x \in A$.
\end{thm}
\begin{pf}
It follows from \eqref{2.4} and \eqref{2.5} that
\begin{eqnarray}\label{2.6}
 \lim_{n\to\infty}2^n\varphi\left(\frac{x}{2}, \frac{y}{2}\right)=0\quad,\quad\lim_{n\to\infty}2^{3n}\psi\left(\frac{a}{2}, \frac{b}{2},\frac{c}{2}\right)=0
\end{eqnarray}
for all $x,y,a,b,c \in V_\bot$. Letting $x=y=a=b=c=0$ in \eqref{2.6} we get $\varphi(0,0)=0,\quad\psi(0,0,0)=0.$
\\Consider the set
$$\Lambda :=\{g:A\to B\quad|\quad g(0)=0\quad g(x)\bot 2g(\frac{x}{2})\quad or \quad 2g(\frac{x}{2})\bot g(x)\}.$$
For every $g,h \in \Lambda$, define,
$$d(g,h):=\inf\{k\in(0,\infty)\quad|\quad\|g(x)-h(x)\|\leq k\varphi(x,x)\quad \forall x\in A\}.$$
Now, put the orthogonally relation $\bot$ on $\Lambda$ as follows: for all $g,h \in \Lambda$
$$h\bot g \Leftrightarrow h(x)\bot g(x)\quad or \quad g(x)\bot h(x)\quad \forall x\in A.$$
It is easy to show $(\Lambda, d, \bot )$ is an O-complete generalized metric space.
\\Now, consider the mapping $T:\Lambda \to \Lambda$ defined by
$$Tg(x):=2g(\frac{x}{2})\quad \forall x\in A.$$
Clearly, T is $\bot$-preserving. As result for all $g,h\in \Lambda$ with $g\bot h$ and $x\in A,$ if $d(g,h)<k,$ then
\begin{equation*}
\begin{split}
&\|g(x)-h(x)\|\leq k\varphi(x,x)\\
&\|2g(\frac{x}{2})-2h(\frac{x}{2})\|\leq2k\varphi(\frac{x}{2},\frac{x}{2})\\
&\|2g(\frac{x}{2})-2h(\frac{x}{2})\|\leq Lk\varphi(\frac{x}{2},\frac{x}{2})\\
& d(Tg,Th)\leq Lk.
\end{split}
\end{equation*}
Hence we see that,
$$d(Tg,Th)\leq Ld(g,h)$$ for all $g,h\in \Lambda$, that is, T is a strictly $\bot$-contractive self mapping of $\Lambda$ with the Lipschitz constant L. The function  $T$ is $\bot$-continuous. In fact, if $\{g_n\}$ is an O-sequence in $\Lambda$ which converges to $g \in \Lambda$, then for given $\varepsilon>0$, there exists $k>0$ with $k<\varepsilon$ and $n \in \mathbb{N}$ such that
 $$\|g_n(x)-g(x)\|\leq k\varphi(x,x)$$
 for all $x\in \mathfrak{A}$ and $n \in \mathbb{N}$.\\
By Theorem \ref{T.1}, there exists a mapping $\Im:A \rightarrow B$ satisfying the following:
\\1) $\Im$ is a fixed point of T, i.e.
  \begin{eqnarray}\label{2.7}
   \Im(x)=2\Im(\frac{x}{2}) \quad \forall x\in A.
\end{eqnarray}
The mapping $\Im$ is a unique fixed point of $T$ in set $\Omega=\{g\in \Lambda : d(f,g)<\infty\}$.\\
$\Im$ is a unique mapping satisfying \eqref{2.7} such that there exists a $k\in (o, \infty)$ satisfying
$$\|f(x)-\Im(x)\|\leq k\varphi (x,x)\quad \forall x\in A.$$
2) $d(T^nf,\Im)\to 0$ as $n\to \infty$. By Theorem \ref{T.1} there exists a fixed point $\Im$ of $T$ such that
\begin{eqnarray}\label{2.8}
\Im(x)=\lim_{n\to \infty}2^nf(\frac{x}{2^n})\quad \forall\quad x\in A
\end{eqnarray}
On the other hand, it follows from \eqref{2.2}, \eqref{2.4} and \eqref{2.6} that
\begin{equation*}
\begin{split}
\|\Delta_\rho \Im(x,y)\|&=\lim_{n\to \infty}2^n\|\Delta_\rho f(\frac{x}{2^n},\frac{y}{2^n})\|\\
&\leq \lim_{n\to \infty}2^n\varphi(\frac{x}{2^n},\frac{y}{2^n})\\
&=0.
\end{split}
\end{equation*}
So $\Delta_\rho \Im(x,y)=0$ for all $x,y \in V_\bot.$ By lemma \ref{l2.1} and [\cite{psl},Theorem 2.2]  $\Delta_\rho \Im(x,y)$ is unique additive. By the result in \cite{v4}, $\Im$ is an orthogonal mapping and so it follows from the definition of $\Im$, \eqref{2.3}, \eqref{2.5} and \eqref{2.6} that we have
\begin{equation*}
\begin{split}
\|\Re \Im([[a,b],c])\|&=\lim_{n\to \infty}2^{3n}\|\Re f([[\frac{a}{2^n},\frac{b}{2^n}],\frac{c}{2^n}])\|\\
&\leq \lim_{n\to \infty}2^{3n}\psi(\frac{a}{2^n},\frac{b}{2^n},\frac{c}{2^n})\\
&=0
\end{split}
\end{equation*}
for all $a,b,c\in V_\bot$.\\
3) $d(f,\Im)\leq\frac{1}{1-L}d(f,Tf)$, which implies that
$$d(f,\Im)\leq\frac{1}{1-L}d(f,Tf)\leq\frac{1}{2(1-L)}\varphi(x,x) \quad \forall x\in A.$$
This completes the proof.
\end{pf}
Now, corollary \ref{c2.3}, shows the Hyers-Ulam-Rassias stability of the orthogonally 3-Lie homomorphism additive $\rho$-Jensen equation \eqref{kv2}.
\begin{cor}\label{c2.3}
Let $s \neq1$ and $\theta$ be nonnegative real numbers, suppose $f :A \rightarrow B$ be is a  mapping  such that
\begin{eqnarray}\label{2.9}
 \|\Delta_\rho f(x,y)\|
\le
 \theta ( \|x\|^s + \|y\|^s )
\end{eqnarray}
and
\begin{eqnarray}\label{2.10}
\|\Re f([[a,b],c])\| \le \theta ( \|a\|^s \cdot \|b\|^s\cdot\|c\|^s )
\end{eqnarray}
for all $x, y, a, b, c\in V_\bot$.
Then there exists a unique 3-Lie homomorphism $\Im : A\rightarrow B$
such that
\begin{eqnarray}\label{2.11}
\|f(x)- A(x) \|  \le  \frac{2 \theta }{2^r -2}  \|x \|^s,\quad\quad\quad(for\quad s<1)
\end{eqnarray}
and
\begin{eqnarray}\label{2.12}
\|f(x)- A(x) \|  \le  \frac{2 \theta }{2 -2^s}  \|x \|^s,\quad\quad\quad(for\quad s>1)
\end{eqnarray}  for all $x \in A$.
\end{cor}

\begin{pf}
The proof follows from Theorem \ref{t2.2} by taking
$$\varphi (x,y)= \theta(\|x\|^r+\|y\|^r)$$
$$\psi(a,b,c)=\theta ( \|a\|^s \cdot \|b\|^s\cdot\|c\|^s$$
for all $x,y,a,b,c\in V_\bot$. Then we can choose  $L= 2^{s-1}$ in \eqref{2.11}, $L= 2^{1-s}$ in \eqref{2.12} and we get the desired result.
\end{pf}
Now, we will prove the hyperstability of orthogonally 3-Lie homomorphism in 3-Lie algebras in the following Theorem.
\begin{thm}\label{t2.4}
Let $f:A\to B$ be a mapping and $\varphi:A^5\to[0, \infty)$ be a function such that
\begin{eqnarray}\label{2.13}
\|\Delta_\rho f(x,y)+ \Re(a,b,c)\|\leq \varphi(0,y,a,b,c)
\end{eqnarray}
for all $x, y, a, b, c \in V_\bot.$ If there exists a constant $0<L<1$ such that
\begin{eqnarray}\label{2.14}
\varphi(0,\frac{y}{2},\frac{a}{2},\frac{b}{2},\frac{c}{2})\leq \frac{L}{2}\varphi(0,y,a,b,c)
\end{eqnarray}
for all $x, y, a, b, c \in V_\bot,$ then f is an orthogonally 3-Lie homomorphism.
\end{thm}
\begin{pf}
Letting $y=a=b=c=0$ in \eqref{2.13} we have
\begin{eqnarray}\label{2.15}
  \left\| f(x) - \frac{1}{2} f(2x) \right\|  \le   \varphi(0,0,0,0,0)
\end{eqnarray}
for all $x \in A$. On the other hand, it follows from \eqref{2.14} that
\begin{eqnarray}\label{2.16}
 \lim_{n\to\infty}2^n\varphi\left(0, \frac{y}{2},\frac{a}{2},\frac{b}{2},\frac{c}{2}\right)=0
\end{eqnarray}
for all $ y, a, b, c \in V_\bot.$ But $\varphi(0,0,0,0,0)=0$, so by \eqref{2.13}, $f(x)=\frac{1}{2} f(2x)$ and then for $n\in \mathbb{N}$ and $x\in A$ we get
\begin{eqnarray}\label{2.17}
f(x)=\frac{1}{2^n} f(2^nx).
\end{eqnarray}
From \eqref{2.13} and \eqref{2.17} we have
\begin{equation}\label{2.18}
\begin{split}
\|\Delta_\rho f(x,y)\|&= \frac{1}{2^n}\|\Delta_\rho f(2^nx,2^ny)\|\\
&\leq \frac{1}{2^n}\varphi\left(0, 2^ny,0,0,0\right)\\
&=0
\end{split}
\end{equation}
for all $x,y \in V_\bot$. So, letting $n\to \infty$ in \eqref{2.18} and using \eqref{2.16}, we have $\|\Delta_\rho f(x,y)\|=0$ for all $x,y \in V_\bot$.\\
On the other hand, we have
\begin{equation}\label{2.19}
\begin{split}
\|\Re f(a,b,c)\|&= \frac{1}{2^n}\|\Re f(2^na,2^nb,2^nc)\|\\
&\leq \frac{1}{2^n}\varphi\left(0, 0,2^na,2^nb,2^nc\right)\\
&=0
\end{split}
\end{equation}
for all $a,b,c \in V_\bot$. Hence, by letting $n\to \infty$  in \eqref{2.19} and using \eqref{2.16}, we have
$$\Re f(a,b,c)=0\quad\quad \forall a,b,c \in V_\bot.$$
Therefore, f is orthogonally 3-Lie homomorphism.
\end{pf}
\begin{cor}\label{c2.6}
Let $\theta$ and $s \neq1$ be nonnegative real numbers. Suppose $f :A \rightarrow B$ is a  mapping  such that
\begin{eqnarray}\label{2.20}
 \|\Delta_\rho f(x,y)+\Re f(a,b,c)\|\le\theta (  \|y\|^s+\|a\|^s +\|b\|^s +\|b\|^s)
\end{eqnarray}
for all $x, y,a,b,c\in V_\bot$. Then f is an orthogonally 3-Lie homomorphism.
\end{cor}

{\footnotesize
\begin{thebibliography}{99}
\bibitem {ul60} S. M. Ulam,  \newblock{\it A Collection of the Mathematical Problems}, \newblock  Interscience Publ. New York, 1960.
\bibitem {hy41} D. H.Hyers,\newblock{\it On  the stability of the linear functional equation}, \newblock Proc. Natl. Acad. Sci. U.S.A., {\bf 27} (1941), 222--224.
\bibitem {ra78} Th. M. Rassias, \newblock{\it  On  the stability of the linear mapping in Banach spaces}, \newblock Proc. Amer. Math. Soc. {\bf 72} (1978), 297--300.
\bibitem{ga94} P. G\v avruta, \newblock{\it A generalization of the Hyers-Ulam-Rassias stability of approximately additive mappings}, \newblock  J. Math. Anal. Appl., {\bf 184} (1994), 431--43.
\bibitem {ir} G.  Isac, Th. M. Rassias, {\it   Stability of $\psi$-additive mappings: Appications to nonlinear analysis},
Internat. J. Math. Math. Sci., {\bf 19} (1996), 219--228.
\bibitem{24}M. Eshaghi Gordji, V. Keshavarz, J. R. Lee, D. Y. Shin  and C. Park, {\it Approximate ternary Jordan ring homomorphisms in ternary Banach Algebras},  J. Comput. Anal. Appl., {\bf 22} (2017), 402--408.
\bibitem{ii2} I.-I. EL-Fassi, S. Kabbaj and A. Charifi, {\it Hyperstability of Cauchy-Jensen functional equations}, Indagationes Mathematicae., {\bf 27} (2016), 855--867.
\bibitem{plx} C. Park, J. Rye Lee and X. Zhang, {\it Additive s-functional inequality and hom-derivations in Banach algebras},  \newblock J. Fixed Point Theory Appl.,  (2019) 14 pages.
\bibitem {mt} C. Park and M. Th. Rassias, {\it Additive functional equations and partial multipliers in $C^*$- algebras},
Revista de la Real Academia de Ciencias Exactas, Serie  A. Matematicas, {\bf 113} (2019), 2261--2275.
\bibitem {mt1} S.-M. Jung and C. Mortici  and M. Th. Rassias , {\it On a functional equation of trigonometric type},
Appled Mathematics and Computation , {\bf 252} (2015), 294--303.
\bibitem {mt2} C. Mortici, S.-M. Jung and M. Th. Rassias, {\it On the stability of a functional equation associated with the Fibonacci numbers},
Abstract and Applied Analysis, Volume 2014 (2014), Article ID 546046, 6 pages.
\bibitem {mt3} S.-M. Jung and M. Th. Rassias, {\it A linear functional equation of third order associated to the Fibonacci numbers},
Abstract and Applied Analysis, Volume 2014 (2014), Article ID 137468, 7 pages.
\bibitem {mt4} S.-M. Jung, D. Popa, and M. Th. Rassias, {\it On the stability of the linear functional equation in a single variable on complete metric groups}, Journal of Global Optimization , {\bf 59} (2014), 165 – 171.
\bibitem {mt5} Y.-H. Lee, S.-M. Jung, and M. Th. Rassias, {\it On an n-dimensional mixed type additive and quadratic functional equation}, Applied Mathematics and Computation, {\bf 228} (2014), 13--16.
\bibitem {Alonso1} J. Alonso, C. Ben\'{\i}tez, {\it Carlos orthogonality in normed linear spaces: A survey. II. Relations between main
orthogonalities}, Extracta Math., {\bf 4} (1989), 121--131.
\bibitem {Alonso2} J. Alonso, C. Ben\'{\i}tez, {\it Orthogonality in normed linear spaces: A survey. I. Main properties}, Extracta Math., {\bf 3} (1988), 1--15.
\bibitem{v3} M. Eshaghi Gordji, M. Ramezani, M. De La Sen, Y.J. Cho, {\it On orthogonal sets and Banach fixed point theorem}, Fixed Point Theory,  {\bf 18}, (2017)  569–-578.
\bibitem {RE} M. Rabbani, M. Eshaghi, {\it Introducing of an orthogonally relation for stability of ternary cubic homomorphisms and derivations on $C^*$-ternary algebras}, Filomat, {\bf 32} (2018), 1439–-1445.
\bibitem{m.r1} M. Ramezani, {\it Orthogonal metric space and convex contractions}, {Int. J. Nonlinear Anal. Appl.}, {\bf 6} (2015), 127--132.
\bibitem{m.r} M. Ramezani and H. Baghani, {\it Contractive gauge functions in strongly orthogonal metric spaces}, {Int. J. Nonlinear Anal. Appl.}, {\bf 8} (2017), 23--28.
\bibitem {SDC} T. Senapati, L. K. Dey, B. Damjanovic and A. Chanda, {\it New Fixed Point Results In Orthogonal Metric Spaces With An Application}, Kragujevac Journal of Mathematics, {\bf 42} (2018), 505--516.
 \bibitem{v4} A. Bahraini, G. Askari, M. Eshaghi Gordji and R. Gholami, {\it Stability and hyperstability of orthogonally $*$-m-homomorphisms in orthogonally Lie $C^*$-algebras: a fixed point approach}, J. Fixed Point Theory Appl., (2018) 12 page.
\bibitem{psl}C. Park, D. Y. Shin and J. R. Lee, {\it Fixed points and additive $\rho$-functional equations}, \newblock J. Fixed Point Theory Appl., {\bf 18} (2016), 569–-586.
\end{thebibliography}
}
\end{document}